\documentclass[11pt]{article}
\usepackage{amsfonts}
\usepackage{amssymb}
\usepackage{enumerate}

\newtheorem{exam}{Example}[section]
\newtheorem{num}  {}

\newcommand{\rmap}{\longrightarrow}
\newcommand{\Boxe}{\raisebox{.8ex}{\framebox}}

\newcommand{\xx}{\ensuremath{\mathcal{X}}}

\newcommand{\A}{\ensuremath{\mathcal{A}}}

\newcommand{\ps}{{\raise 1pt\hbox{\tiny (}}}

\newcommand{\pss}{{\raise 1pt\hbox{\tiny [}}}
\newcommand{\pdd}{{\raise 1pt\hbox{\tiny ]}}}
\newcommand{\pd}{{\raise 1pt\hbox{\tiny )}}}

\newcommand{\bs}{{\raise 1pt\hbox{\tiny [}}}
\newcommand{\bd}{{\raise 1pt\hbox{\tiny ]}}}

\def\cross{\mathinner{\mathrel{\raise0.8pt\hbox{$\scriptstyle>$}}
                 \joinrel\mathrel\triangleleft}}

\def\compose{{\raise 1pt\hbox{$\scriptscriptstyle\circ$}}}
\def\dcross{{\raise 0.5pt\hbox{$\scriptscriptstyle\boxtime$}}}

\advance\textheight by \topskip
\usepackage[all]{xy}
\usepackage{epsf}
\usepackage{amsfonts}
\usepackage{amssymb}
\usepackage{amsmath}

\begin{document}

\title{Chern characters via connections up to homotopy \thanks{Research supported by NWO}}
\author {Marius Crainic}
\date {Department of Mathematics, Utrecht University, The Netherlands}
\pagestyle{myheadings}
\maketitle
%\begin{abstract}
%In this note we remark that Chern characters can be computed using curvatures of ``(super-)connections up to homotopy''.
%As an application, we prove the vanishing theorem for Lie algebroids.
%\end{abstract}
%********

%\tableofcontents

\begin{num}{\bf Introduction: }\rm The aim of this note is to point out that Chern characters can be computed using
curvatures of ``connections up to homotopy'', and to present an application to the vanishing theorem for Lie algebroids.\\
\hspace*{.3in}Classically, Chern characters are computed with the help of
a connection and its curvature. However, one often has to relax the notion of connection so that one gains more freedom
in representing these characteristic classes by differential forms. A well known example is Quillen's notion of super-connection \cite{Chern}.
Here we remark that one can weaken the notion of (super-)connections
even further, to what we call ``up to homotopy''. \\
\hspace*{.3in}Our interest on this type of connections comes from the theory of characteristic classes of Lie algebroids \cite{Cra, Fer}
(hence, in particular of Poisson manifolds \cite{Fer2}). From our point of view, the intrinsic characteristic classes
are secondary classes which arise from a vanishing result:
the Chern classes of the adjoint representation vanish (compare to Bott's approach to characteristic classes for foliations).
We have sketched a proof of this in \cite{Cra} for a particular class of Lie algebroids (the so called regular ones).
The problem is that the adjoint representation is a representation up to homotopy only \cite{ELW}. For the general setting, 
we have to show that Chern classes can be computed using connections up to homotopy. Since we believe that this result might be 
of larger interest, we have chosen to present it at the level vector bundles over manifolds. In \cite{Crai} we
will describe the secondary characteristic classes which arise in the flat case.
\end{num}

\hspace*{.1in}All objects in this note should be viewed in the super (i.e. $\mathbb{Z}_2$-graded) setting. E.g. we work with super-vector spaces, super-algebras, super-commutators (i.e. $[a, b]= ab - (-1)^{|a|\cdot |b|} ba$), super-traces 
(i.e. $Tr_s(A)= Tr(A_{0\, 0})- Tr(A_{1\, 1})$ for any
endomorphism $A\in \text{End}(V)$ of a super-vector space $V= V^0\oplus V^1$).\\
\hspace*{.3in}Let $(E, \partial)$ be a super-complex of vector bundles over a manifold $M$,
\begin{eqnarray}\label{complex}
(E, \partial): \ \  \xymatrix{  E^0\ \ar@<-1ex>[r]_-{\partial} & \ \ E^1 \ar@<-1ex>[l]_-{\partial}\ \ \ \ . }
\end{eqnarray}
Such an object can be viewed as an element in the $K$-theory of $M$ (the formal differences $E^0- E^1$). Accordingly,
the Chern character of $E$ is
\[ Ch(E)= Ch(E^0)- Ch(E^1) \in H^{*}(M) \ .\]

\begin{num}{\bf Connections: }\rm A {\it connection} on $E$ is a linear map
\begin{equation}\label{connection}
 \xx(M)\otimes \Gamma E \rmap \Gamma E, \ \ (X, s)\mapsto \nabla_{X}(s)
\end{equation}
with the following properties:
\begin{enumerate} [(i)]
\item $\nabla_X$ is even (i.e. preserves the degrees), and $\nabla_X\partial= \partial\nabla_X$\ ,
\item $\nabla_{X}(fs)= f\nabla_{X}(s) + X(f) s$\ ,
\item $\nabla_{fX}(s)= f\nabla_{X}(s)$\ ,
\end{enumerate}
for all $X\in \xx(M)$, $s\in \Gamma E$, and $f\in C^{\infty}(M)$.
\end{num}

\begin{num}{\bf Connections up to homotopy: }\rm A {\it connection up homotopy} on $E$ is a
local operator (\ref{connection}) which satisfies the properties (i) and (ii)
above, and satisfies (iii) up to homotopy only. In other words we require
\begin{equation}\label{iiihom}
\nabla_{fX}= f\nabla_{X} + [H_{\nabla}(f, X), \partial] \ ,
\end{equation}
where $H_{\nabla}(f, X)\in \Gamma \text{End}(E)$ are local operators of odd degree, linear on $X$ and $f$.
\end{num}

\begin{num}{\bf Nonlinear forms: }\rm Many of the basic operations on the space $\A(M; E)$ of
$E$-valued differential forms
\begin{equation}\label{omega}
\omega: \underbrace{\xx(M)\times \ldots \times \xx(M)}_{n} \rmap \Gamma(E)
\end{equation}
hold without any $C^{\infty}(M)$-linearity assumption on $\omega$.
We recall these basic operations. So, let us denote by
$\A_{\text{nl}}^{n}(M; E)$ the space of all antisymmetric
($\mathbb{R}$-multilinear) maps (\ref{omega}). The familiar
formula
\begin{eqnarray}\label{product}
& (\omega \eta)(X_1, \ldots , X_{n+m}) = & \nonumber \\
& = \sum_{\sigma\in S(n, m)} \text{sgn}(\sigma) \omega(X_{\sigma(1)}, \ldots ,
X_{\sigma(n)}) \eta(X_{\sigma(n+1)}, \ldots , X_{\sigma(n+m)}) &
\end{eqnarray}
(where $S(n, m)$ stands for $(n, m)$-shuffles) extends the usual
product of forms to a product on $\A_{\text{nl}}^{n}(M)$. The same
formula defines a left action of $\A_{\text{nl}}(M)$ on $\A_{\text{nl}}(M;
E)$.\\ 
\hspace*{.3in}Any local operator $\nabla$ satisfying (i)
and (ii) above can be viewed as a map $\A_{\text{nl}}^{0}(M; E)\rmap
\A_{\text{nl}}^{1}(M; E)$, and it has (by the classical arguments) a
unique extension to an odd operator
\begin{equation}\label{operator}
\nabla: \A_{\text{nl}}(M; E)\rmap \A_{\text{nl}}(M; E)
\end{equation}
which satisfies the Leibniz rule. Explicitly,
\begin{eqnarray}\label{differential}
\nabla(\omega)(X_1, \ldots , X_{n+1}) & = & \sum_{i<j}
(-1)^{i+j}\omega([X_i, X_j], X_1, \ldots , \hat{X_i}, \ldots ,
\hat{X_j}, \ldots X_{n+1})) \nonumber \\
 & + & \sum_{i=1}^{n+1}(-1)^{i+1}
\nabla_{X_i}\omega(X_1, \ldots, \hat{X_i}, \ldots , X_{p+1}) .
\end{eqnarray}
In the particular case where $E$ is trivial (and $\nabla_X$ is the
Lie derivate along $X$), this gives an operator $d$ on
$\A_{\text{nl}}(M)$ which extends the classical De Rham operator $d$ on
differential forms.\\ 
\hspace*{.3in}Note that $\A_{\text{nl}}(M; E)$ is
canonically isomorphic to
$\A_{\text{nl}}(M)\otimes_{C^{\infty}(M)}\Gamma(E)$. In particular,
$\A_{\text{nl}}(M; \text{End}(E))$ has a canonical product. This product is
given by the same formula (\ref{product}) except for a minus sign
(due to the definition of the tensor product of super-algebras
\cite{Chern}) which appears when $m$ and the $E$-degree of
$\omega$ are odd.\\ 
\hspace*{.3in}For $\nabla$ as above, the operators $[\nabla_X, T]$ acting on $\Gamma(E)$ are
$C^{\infty}(M)$-linear for any $T\in \Gamma \text{End}(E)$, and the
correspondence $(X, T)\mapsto [\nabla_X, T]$ defines a similar
operator $\tilde{\nabla}$ on $\Gamma \text{End}(E)$. Its extension to
$\A_{\text{nl}}(M; \text{End}(E))$ is denoted by
\[ d_{\nabla}: \A_{\text{nl}}(M; \text{End}(E)) \rmap \A_{\text{nl}}(M; \text{End}(E)) \]
Identifying the elements of $\A_{\text{nl}}(M; \text{End}(E))$ with the induced multiplication operators
(acting on $\A_{\text{nl}}(M; E)$), $d_{\nabla}(-)$ coincides with the (graded) commutator $[\nabla, -]$ with (\ref{operator}). \\
\hspace*{.3in}The square of (\ref{operator}) is the product by an element
$k_{\nabla}\in \A_{\text{nl}}(M; \text{End}(E))$, {\it the curvature of $\nabla$}, which
is given by the usual formula
\begin{equation}\label{curvature}
k_{\nabla}(X, Y)= [\nabla_{X}, \nabla_{Y}]- \nabla_{[X, Y]}: \Gamma E\rmap \Gamma E
\end{equation}
\end{num}

\hspace*{-.2in}{\bf Lemma 1. }{\it For any local operator $\nabla$ satisfying (i) and (ii) above,
\begin{enumerate} [(i)]
\item $\A_{{\rm nl}}(M;\, {\rm End}(E))$ is a ($\mathbb{Z}_{2}$-graded) algebra endowed with an odd operator $d_{\nabla}$
which satisfies the Leibniz rule;
\item $d_{\nabla}^{2}$ is the commutator by $k_{\nabla}$, and $d_{\nabla}(k_{\nabla})= 0$;
\item The super-trace defines a map
\begin{equation}\label{supertrace}
Tr_{s}: (\A_{{\rm nl}}(M;\, {\rm End}(E)), d_{\nabla}) \rmap (\A_{{\rm nl}}(M), d)
\end{equation}
with the property that $Tr_{s}d_{\nabla}= d Tr_{s}$.
\item For any $p\geq 0$, $Tr_{s}(k_{\nabla}^{p})\in \A_{{\rm nl}}^{2p}(M)$ is closed.
\end{enumerate}}

{\it Proof: }As in the classical case \cite{Chern}, the last part follows from (i)-(iii). We still have to prove (iii). Since by adding new vector bundles (and extend $\nabla$ with the help of
any connections) we can make $E^0$ and $E^1$ trivial, we may assume that $E$ is trivial as a vector bundle.
In this case the assertion follows easily from the fact that
$d_{\nabla}= d+ [\theta, -]$ for some $\theta\in \A_{\text{nl}}^{1}(M; \text{End}(E))$.  \ \ $\Boxe$\\

%\hspace*{.1in}As usual (see e.g. \cite{Chern}), given two operators $\nabla^0$, $\nabla^1$ as above, 
%$Tr_{s}(k_{\nabla{0}}^{p})$ and $Tr_{s}(k_{\nabla_{1}}^{p})$ differ by a boundary.
%Indeed, we can form the affine combination $\nabla_{t}= (1-t)\nabla_0+ t\nabla_1$ acting 
% on the pull-back of $E$ to $M\times I$ ($I$ is the unit interval),
%and then the Chern-Simons type forms $cs_{p}(\nabla_0, \nabla_1)= \int_{0}^{1} Tr_{s}(k_{\nabla_t}^{p})$
%satisfy
%\begin{equation}\label{CSim}
%Tr_{s}(k_{\nabla{0}}^{p}) - Tr_{s}(k_{\nabla_{1}}^{p})= d cs_{p}(\nabla_0, \nabla_1)\ .
%\end{equation}

\begin{num}{\bf Forms up to homotopy: }\rm There are various cases where the Chern-type
elements $Tr_{s}(k_{\nabla}^{p})$
of the previous lemma are true forms on $M$. This happens for instance when $\nabla$ is a
connection up to homotopy.
To see this, we consider $\text{End}(E)$-valued
{\it forms up to homotopy}, i.e. elements $\omega\in \A_{\text{nl}}(M; \text{End}(E))$ which commute with
$\partial$ and are ``$C^{\infty}(M)$- linear up to homotopy''.
In other words, we require that
\[ f\omega(X_1, \ldots , X_n)- \omega(fX_1, \ldots , X_n)=
[H_{\omega}(f, X_1, \ldots , X_n), \partial ] \] 
for some operator
\[ H_{\omega}(f, X_1, \ldots , X_n)\in \Gamma \text{End}(E)\] 
depending linearly on $f\in C^{\infty}(M)$ and on the vector fields $X_i$.
We denote by $\A_{\partial}^{n}(M; \text{End}(E))$ the space of such
$\omega$'s.
\end{num}

\hspace*{-.2in}{\bf Lemma $2$. }{\it For any connection up to homotopy $\nabla$,
\begin{enumerate} [(i)]
\item $k_{\nabla}\in \A_{\partial}^{2}(M; \, {\rm End}(E))$ ;
\item $\A_{\partial}^{n}(M; \, {\rm End}(E))$ is a subalgebra of $\A_{{\rm nl}}(M; \, {\rm End}(E))$, which
is preserved by $d_{\nabla}$;
\item The trace map (\ref{supertrace}) restricts to
\[ Tr_s: (\A_{\partial}(M; \, {\rm End}(E)), d_{\nabla}) \rmap (\A(M), d) \ .\]
\end{enumerate}}

\hspace*{.1in}The proof is a simple (and standard) computation.
And now the conclusion:\\

\hspace*{-.2in}{\bf Theorem. }\label{theorem}{\it If $\nabla$ is a connection up to homotopy
on $(E, \partial)$, and $k= k_{\nabla}$ is its curvature, then
\begin{equation}\label{Chcl}  Tr_{s}(k^{p}) \in \A^{2p}(M)
\end{equation}
are closed forms whose De Rham cohomology classes are (up to a constant) the components of
the Chern character $Ch(E)$.}\\

{\it Proof: }We still have to show that (\ref{Chcl}) induce the components of the Chern
character. This is clear (by the definition of the Chern character) if $\nabla$ is a connection on $(E, \partial)$. 
In general, we choose such a connection $\widetilde{\nabla}$. To see that such $\widetilde{\nabla}$ exists,
we can locally define $\widetilde{\nabla}_{f \frac{\partial}{\partial x_i}}= f \nabla_{\frac{\partial}{\partial x_i}}$
and then use a partition of unity argument. Hence it suffices to show that $Tr_{s}(k_{\nabla}^{p})$ and $Tr_{s}(k_{\widetilde{\nabla}}^{p})$
differ by the differential of a (true) differential form on $M$. For this we form the affine combination $\nabla_{t}= (1-t)\nabla+ t\widetilde{\nabla}$
which is a connection up to homotopy on the pull-back of $(E, \partial)$ to $M\times I$ ($I$ is the unit interval). The 
Chern-Simons type forms $cs_{p}(\nabla, \widetilde{\nabla})= \int_{0}^{1} Tr_{s}(k_{\nabla_t}^{p})$ will satisfy the desired 
equation $Tr_{s}(k_{\nabla}^{p}) - Tr_{s}(k_{\widetilde{\nabla}}^{p})= d cs_{p}(\nabla, \widetilde{\nabla})$.\ \ \ $\Boxe$\\

\begin{num}{\bf Application to Lie algebroids: }\rm Recall \cite{McK} that a
{\it Lie algebroid} over $M$ is a triple
\[ (\mathfrak{g}, \, [\cdot \, , \cdot ] \, , \rho) \]
consisting of a vector bundle $\mathfrak{g}$ over $M$, a Lie bracket
$[ \cdot \, , \cdot ]$ on the space $\Gamma(\mathfrak{g})$, and a
morphism of vector bundles $\rho: \mathfrak{g}\rmap TM$ ({\it the anchor}
of $\mathfrak{g}$), so that $[X, fY] = f[X,Y]+ \rho(X)(f) \cdot Y$ for all
$X, Y \in \Gamma(\mathfrak{g})$ and $f \in C^{\infty}(M)$. Basic examples are Lie algebras
(when $M$ is a point), the tangent bundle (when $\mathfrak{g}= TM$ and $\rho$ is the identity), and foliations (when $\rho$ is the inclusion
of the involutive bundle of vectors tangent to the foliation). 
Formula (\ref{differential}) defines a ``$\mathfrak{g}$-De Rham differential''
 on the space $C^*(\mathfrak{g})= \Gamma(\Lambda \mathfrak{g}^*)$, and the resulting
cohomology is well known  under the name of {\it the cohomology of the Lie algebroid} $\mathfrak{g}$, denoted by $H^*(\mathfrak{g})$ \cite{McK}.
It generalizes
Lie algebra cohomology and De Rham cohomology. It also relates to the last one by composition with
the anchor map:
\begin{equation}\label{anchor}
\rho_{*}: H^{*}(M)\rmap H^{*}(\mathfrak{g})\ .
\end{equation}
\end{num}
\hspace*{.3in}The formal difference $TM- \mathfrak{g}$ plays the role of the ``normal bundle'' of $\mathfrak{g}$
(and is minus the ``adjoint representation'' \cite{ELW}). The following is at the origin of the appearance of 
secondary characteristic classes for Lie algebroids \cite{Cra, Fer}. \\

\hspace*{-.2in}{\bf Corollary. }{\it (Vanishing theorem) For any Lie algebroid $\mathfrak{g}$,
(\ref{anchor}) kills the Chern character of the difference bundle $\mathfrak{g}- TM$.}\\

{\it Proof: } There is an obvious notion of $\mathfrak{g}$-connection on a vector bundle $E$:
a map $\Gamma(\mathfrak{g})\times \Gamma E\rmap \Gamma E$ satisfying the usual identities. 
Exactly as in the classical case, $\mathfrak{g}$-connections define cohomology classes
$Ch_{\mathfrak{g}}(E)\in H(\mathfrak{g})$, which do not depend on $\nabla$.
Since any (classical) connection on $E$ induces an obvious $\mathfrak{g}$ connection, we see
that $Ch_{\mathfrak{g}}(E)= \rho_{*}Ch(E)$. On the other hand, all our arguments have an obvious
$\mathfrak{g}$-version (just replace $\xx(M)$ by $\Gamma(\mathfrak{g})$). Hence it suffices
to point out \cite{ELW} that the super-complex
\begin{equation}\label{adjrep}
\xymatrix{  \mathfrak{g}\ \ar@<-1ex>[r]_-{\rho} & \ \ TM \ar@<-1ex>[l]_-{0}\ \ \ \ . }
\end{equation}
is endowed with a canonical $\mathfrak{g}$-connection up to homotopy which is flat:
\begin{eqnarray}
\nabla_{X}(Y)= [X, Y]\ , & \nabla_{X}(V)= [\rho(X), Y] \ ,  \nonumber \\
H(f, X)(Y)= 0\ , & H(f, X)(V)= V(f) X \ .  \nonumber
\end{eqnarray}
($X, Y\in \Gamma\mathfrak{g}$, $V\in \xx(M)$, $f\in C^{\infty}(M)$). \ \ $\Boxe$\\

\begin{num}{\bf Variations: }\rm
\begin{enumerate}[(i)]
\item  Exactly the same discussion applies to  super-connections \cite{Chern} as well (note that our exposition is designed
for that). Any operator $\nabla: \A_{\text{nl}}(M; E)\rmap \A_{\text{nl}}(M; E)$ of degree one (with respect to the even/odd grading) which satisfies the Leibniz rule has a decomposition
\[ \nabla= \nabla_{[0]} + \nabla_{[1]}+ \nabla_{[2]} + \ldots\ ,\]
where $\nabla_{[1]}: \xx(M)\times \Gamma E\rmap \Gamma E$, and $\nabla_{[i]}$ are multiplication by elements $\omega_{[i]}\in \A_{\text{nl}}^{i}(M; \text{End}(E))$. We say that $\nabla$ is a {\it super-connection up to homotopy} if $\nabla_{[1]}$ is a connection
up to homotopy and $\nabla_{[i]}$ are multiplication by elements $\omega_{[i]}\in \A_{\partial}^{i}(M; \text{End}(E))$.
There is an obvious version of our discussion for such $\nabla$'s (and a non-linear version of \cite{Chern}).
\item If (\ref{complex}) is contractible over $A\subset M$ (i.e. defines an element in the relative $K$-theory groups),
then the Chern-character forms given by our theorem vanish when restricted to $A$, hence define relative cohomology classes.
\item Note that a consequence of (\ref{iiihom}) is that
\[ [fH^1(g, X)- H^{1}(fg, X)+ H^{1}(f, gX), \partial ]= 0  \ ,\]
for all $f, g\in C^{\infty}(M)$, $X\in \xx(M)$. Hence it seems natural to require that
\begin{equation}\label{higher}
fH^1(g, X)- H^{1}(fg, X)+ H^{1}(f, gX)= [ H^2(f, g, X), \partial ]
\end{equation}
where $H^{2}(f, g, X)\in \Gamma \text{End}(E)$. Similarly, it seems natural to require the existence
of higher $H^{k}$'s satisfying the higher versions of (\ref{higher}). A similar remark applies also to the definition of forms up to homotopy, and induces a subalgebra of $\A_{\partial}(M; \text{End}(E))$ with similar properties.
\end{enumerate}
\end{num}

Marius Crainic,\\
\hspace*{.2in}Utrecht University, Department of Mathematics,\\
\hspace*{.2in}P.O.Box:80.010,3508 TA Utrecht, The Netherlands,\\ 
\hspace*{.2in}e-mail: crainic@math.ruu.nl\\
\hspace*{.2in}home-page: http://www.math.uu.nl/people/crainic/

\end{document}